\documentclass[12pt,amscd]{amsart}
\usepackage[all]{xy}
\usepackage{graphicx}
\usepackage{amsmath,amsxtra,amssymb,latexsym, amscd,amsthm}
\usepackage{indentfirst}
\usepackage[mathscr]{eucal}

\newtheorem{thm}{Theorem}[section]

\newtheorem{lem}[thm]{Lemma}
\newtheorem{prop}[thm]{Proposition}
\theoremstyle{definition}
\newtheorem{defn}[thm]{Definition}
\newtheorem{exm}[thm]{Example}
\newtheorem{rem}[thm]{Remark}

\numberwithin{equation}{section}

\DeclareMathOperator{\depth}{depth}
\DeclareMathOperator{\ass}{Ass}
\DeclareMathOperator{\ann}{ann}
\DeclareMathOperator{\Ext}{Ext}

\DeclareMathOperator{\supp}{supp}
\DeclareMathOperator{\rank}{rank}

\DeclareMathOperator{\dstab}{dstab}
\def\dstabIt{\overline{\operatorname{dstab}}}
\DeclareMathOperator{\astab}{astab}
\def\astabIt{\overline{\operatorname{astab}}}
\DeclareMathOperator{\conv}{conv}
\DeclareMathOperator{\Min}{Min}
\DeclareMathOperator{\height}{ht}
\DeclareMathOperator{\bight}{bight}
\DeclareMathOperator{\link}{lk}

\def\Nset{\mathbb {N}}
\def\Zset{\mathbb {Z}}

\def\Rset{\mathbb {R}}
\def\Ga {CS_{\albf}}
\def\Da {\Delta_{\albf}}

\def\Fcal{\mathcal F}
\def\Rcal{\mathcal R}
\def\Scal{\mathcal S}
\def\Ecal{\mathcal E}
\def\Hcal{\mathcal H}
\def\Bcal{\mathcal B}
\def\Jcal{\mathcal J}

\def\albf {{\boldsymbol{\alpha}}}
\def\bebf {{\boldsymbol{\beta}}}

\def\zbf {\mathbf 0}
\def\ebf {\mathbf e}
\def\xbf {\mathbf x}
\def\Xbf {\mathbf X}

\def\abf {\mathbf a}

\def\mfr {\mathfrak m}

\def\pfr {\mathfrak p}
\def\afr {\mathfrak a}
\def\Fcal{\mathcal F}

\begin{document}

\title[Powers of  Ideals and Combinatorics] {Powers of Monomial Ideals and Combinatorics}
\author{Le Tuan Hoa }
\address{Institute of Mathematics, VAST, 18 Hoang Quoc Viet, 10307 Hanoi, Viet Nam}
\email{lthoa@math.ac.vn}
\subjclass{13D45, 05C90}
\keywords{Associated prime, depth, monomial ideal,  integral closure, simplicial complex,  integer linear programming}
\date{}
\dedicatory{Dedicated to the 80th birthday of Professor L. Bokut.}
\commby{}
\begin{abstract} This is an exposition of some new results on associated primes and the depth  of  different kinds of powers of  monomial ideals in order to show a deep connection between commutative algebra and some objects in combinatorics such as simplicial complexes, integral points in polytopes  and graphs.
\end{abstract}

\maketitle
\section*{Introduction}

The interaction between commutative algebra and combinatorics has a long history. It goes back at least to Macaulay's article \cite{Mac}. Stanley's solution of the so-called Upper Bound Conjecture for spheres gives a new impulse for the study in this direction. Since then, many books devoted to  various topics of this interaction are published, see, e.g.,  \cite{St, BH, Vi,  HH2, MS}. People even talk about the birth of a new area of mathematics called ``Combinatorial Commutative Algebra".

This exposition is based on my talk at the ``Third International Congress in Algebras and Combinatorics (ICAC 2017)'' held Hong Kong. The aim of the workshop is  clear from its title: to know better about the interaction between various areas of mathematics, which include associative algebra, commutative algebra and combinatorics. So, the purpose of this paper is to provide some further  interaction from current research interest. Two basic notions in commutative algebra are concerned here: the associated primes  and the depth of a (graded or local) ring. Associated prime ideals of a ring play a role like prime divisors of a natural number in the Number Theory, while the depth measures how far the ring from being Cohen-Macaulay.

For simplicity, we  work with homogeneous ideals $I$ in a polynomial ring $R= K[X_1,...,X_r]$. On the way to give a counter-example to Conjecture 2.1 in \cite{Ra0}, Brodmann \cite{B1} proves that $\ass(R/I^n)$ becomes stable for $n\gg 0$. This stable set is denoted by $\ass^\infty (R/I)$. Since associated primes are closely related to the depth, almost at the same time Brodmann \cite{B2} proves that $\depth (R/I^n)$ becomes constant for all $n\gg 0$.  This constant is denoted by $\lim_{n\rightarrow \infty }\depth (R/^n)$. 
It is however not known, when the sequences $\{ \ass(R/I^n) \}$ and $\{ \depth(R/I^n)\}$ become stable, and it is  little known about $\ass^\infty (R/I)$ and $\lim_{n\rightarrow \infty }\depth R/^n$. Therefore, it is of great interest to bound the least place $\astab(I)$ (resp. $\dstab(I)$) when the stability of $\ass(R/I^n)$ (resp. $\depth (R/I^n)$) occurs (see Definition \ref{ast} and Definition \ref{dst}), as well as to determine $\ass^\infty (R/I)$ and $\lim_{n\rightarrow \infty }\depth (R/I^n)$. For an arbitrary ideal, these problems are very difficult, because  it is not known how to compute all associated primes of a ring and there is  no effective way to compute the depth. Only few restricted results were obtained in the general case, see \cite{ME, Sh} and \cite[Theorem 2.2]{Mo}.

Luckily, in the case of monomial ideals (i.e. ideals generated by monomials) one can use combinatorics to compute associated primes as well as the depth.   Even it is not a trivial task, it opens up a way to use combinatorics to deal with  these problems. This paper is mainly devoted to bounding $\astab(I)$ and $\dstab(I)$ for monomial ideals. This topic attracts many researchers during  the  last decade. From this study some times one can get surprising relationships between
 seemingly unrelated notions of commutative algebra and combinatorics. For an example, \cite[Theorem 1.2]{TT1} states that  the ring $R/I^n_\Delta$ is Cohen-Macaulay for some fixed $n\ge 3$ if and only if the simplicial complex $\Delta$  is a complete intersection (see Theorem \ref{TTThm2}).  Together with stating some main results we also give some hints for their proofs. Techniques from combinatorics used to obtain results presented in this paper are so broad, that in the most cases we cannot go to the details. We only explain in more details how the existence of integer solutions of  systems of linear constrains related to  bounding $\astab(I)$ and $\dstab(I)$. Besides these two problems, we also list some results on  properties of the sequences $\{\ass(R/I^n) \}$ and $\{\depth (R/I^n) \}$, because they are useful in determining  $\ass^\infty (R/I)$ and $\lim_{n\rightarrow \infty }\depth R/^n$.  Similar problems for integral closures of powers as well as symbolic powers are also considered in this paper. 

We would like to mention that recently there is an intensive research on the so-called Castelnuovo-Mumford regularity of monomial ideals, which also involves a lot of combinatorics. The interested readers can consult the survey paper \cite{BBT}. 

The paper is organized as follows. In Section \ref{Pre} we recall some basic notions and facts from commutative algebra and formulate two main problems considered in this paper. In particular, the above mentioned Brodmann's results are stated here. Section \ref{Ass} is devoted to bounding $\astab(I)$. This section is divided to three subsections: bounds on $\astab(I)$ and $\astabIt(I)$ are presented in the first two subsections. These bounds are huge ones. Good bounds on these invariants for some classes of  monomial ideals are given in the last subsection. The stability of the depth function is presented in Section \ref{Depth}. The first three subsections are devoted to three kind of  powers. The last subsection is concerning with Cohen-Macaulay property of square-free monomials ideals.

\section{Preliminaries} \label{Pre}

Let $R$ be either a Noetherian local ring with maximal ideal $\mfr$ and $K=R/\mfr$, or a standard graded finitely generated $K$-algebra with graded maximal ideal $\mfr$, where $K$ is an infinite field (standard grading means $R= \oplus_{i\ge 0}R_i$ such that $R_0 = K$, $R_iR_j \subseteq R_{i+j}$ for all $i,j\ge 0$   and $R$ is generated by $R_1$ over $K$). A non-zero divisor $x\in R$ is called an $R$-regular element. A sequence of elements $x_1,...,x_s$ of $R$ is called $R$-regular sequence if $x_i$ is an $R/(x_1,...,x_{i-1})$-regular element for $i=1,...,s$, and $R \neq (x_1,...,x_s)$. Then all maximal $R$-regular sequences in $\mfr$ have the same length and this length is called the {\it depth} of $R$, denoted by $\depth (R)$.  Moreover, in the graded case, one can choose a maximal $R$-regular sequence consisting of homogeneous elements. Rees \cite{Re} showes that 
$$\depth(R)= \min\{i|\ \Ext^i_R(K, R) \neq 0\}.$$
One can also define $\depth(R)$ by using  local cohomology:
$$ \depth (R) = \min\{i|\ H^i_{\mfr}( R) \neq 0 \}.$$
The reader can consult the book \cite{BS} for the definition and a detailed algebraic introduction to Grothendieck's local cohomology theory. 

The Krull dimension $\dim (R)$ of $R$ and $\depth (R)$ are two basic invariants of $R$. One has $\depth R \le \dim R$. When the equality holds, $R$ is called a Cohen-Macaulay ring. ``The notion of (local) Cohen-Macaulay ring is a workhorse of Commutative Algebra", see \cite[p. 56]{BH}. This explains the importance of depth.

 It is in general not easy to determine the exact value of $\depth (R)$. Therefore, the following simple result of Brodmann \cite{B2} is of great interest. 

\begin{thm} {\rm (\cite[Theorem 2]{B2})} \label{Brod1} Let $I \subset R$ be a proper ideal, which is assumed to be graded if $R$ is graded. Then

 (i)  $\depth (R/I^n)$ is constant for all $n\gg 0$. 

(ii) Denote the above constant  by $\lim_{n\rightarrow \infty }\depth(R/I^n)$. Let $\Rcal(I) = \oplus_{n\ge 0} I^nt^n$ be the Rees algebra of $I$. Then
$$\lim_{n\rightarrow \infty }\depth(R/I^n) \le \dim (R) - \ell(I),$$
where $\ell (I) = \dim  \Rcal(I)/\mfr \Rcal$ is the analytic spread of $I$.
\end{thm}

In fact, Brodmann's result was formulated for modules.  Brodmann's proof  as well as a new proof by Herzog and Hibi (see \cite[Theorem 1.1]{HH1})  are based on the Noetherian property of the Rees algebra $\Rcal(I)$.  As a corollary of these proofs, one has a similar statement for the so-called integral closures of powers of an ideal. Recall that the {\it integral closure} of an arbitrary ideal $\afr$ of $R$ is the set of elements $x$ in $R$ that satisfy an integral relation
$$x^n + a_1x^{n-1} + \cdots + a_{n-1}x + a_n = 0,$$
where $a_i \in \afr^i$ for $i = 1,\ldots, n$. This is an ideal of $R$ and is denoted by $\overline \afr$. In the local case, assume in addition that $R$ is complete. Then the algebra $\overline{\Rcal}(I) :=  \oplus_{n\ge 0}\overline{ I^n}t^n$ is a module-finite extension of $\Rcal(I)$. So, Brodmann's result implies that $\depth(R/ \overline{ I^n})$ also is constant for all $n\gg 0$.

\begin{defn} \label{dst} For an ideal $I\subset  R$, set
$$\dstab (I) := \min\{ s|\ \depth(R/I^n) = \depth (R/I^s) \ \forall n\ge s \} .$$
In the local case, assume in addition that $R$ is complete. Set
$$\dstabIt (I) := \min\{ s|\ \depth(R/ \overline{ I^n}) = \depth (R/ \overline{ I^s}) \ \forall n\ge s \} .$$
 \end{defn}
 
 The proofs of Brodmann and Herzog-Hibi  give no information on when the functions $\depth (R/I^n)$ and  $\depth(R/ \overline{ I^n})$  become stable. Therefore, the following problem attracts attention of many researchers:
 
 \vskip0.3cm
 \noindent {\bf Problem 1}. {\it Give  upper bounds on $\dstab(I)$ and $\dstabIt(I)$ in terms of other invariants of $R$ and $I$.}
 \vskip0.3cm
 
 Until now there is no approach to solve this problem in the general setting as above. The reason is that there is no effective way to compute depth. Therefore all known nontrivial results until now are dealing with monomial ideals in a polynomial ring.  These results will be summarized in Section \ref{Depth}. Below we describe one of the main tools to be used. 

For the moment, let $R= K[X_1,...,X_r]$ be a polynomial ring with $r$ indeterminate $X_1,..., X_r$. A monomial ideal $I$ of $R$ is an ideal generated by monomials ${\Xbf}^{\albf} := X_1^{\alpha_1}\cdots X_r^{\alpha_r}$, where $\albf = (\alpha_1,...,\alpha_r) \in \Nset^r$. In this case one can effectively describe the local cohomology module $H^i_\mfr(R/I)$, where $\mfr = (X_1,...,X_r)$. Let us recall it here.
  
Since $R/I$ is an $\Nset^r$-graded algebra, $H_{\mfr}^i(R/I)$ is an $\Zset^r$-graded module over $R$, i.e. $H_{\mfr}^i(R/I) = \oplus_{\albf \in \Zset^r} H_{\mfr}^i(R/I)_\albf$, such that $\Xbf^\bebf H_{\mfr}^i(R/I)_\albf \subseteq H_{\mfr}^i(R/I)_{\albf + \bebf}$.  Each $\albf$-component  $H_{\mfr}^i(R/I)_\albf$ can be computed via the reduced simplicial homology. 

 Recall that a simplicial complex  $\Delta$ on the finite set $[r] := \{1,\ldots,r\}$ is a collection of subsets of $[r]$ such that $F\in \Delta $ whenever $F\subseteq F'$ for some $F'\in \Delta $. Notice that we do not impose the condition that $\{i\} \in \Delta $ for all $i\in [r]$.  An element of $\Delta $ is called  a {\it face}. A simplicial complex $\Delta $ is defined by  the set of its {\it facets} (i.e. maximal faces) - denoted by $\Fcal(\Delta)$. In this case we also write $\Delta = \langle \Fcal(\Delta) \rangle$. To each monomial ideal $I$ we can associate a simplicial complex  $\Delta(I)$ defined by 
$$\Delta (I) =\{ \{i_1,...,i_s\} \subseteq [r]|\ X_{i_1}\cdots X_{i_s} \not\in \sqrt I\}.$$
Thus $ \Delta (I)$ is defined upto the radical $\sqrt I$ of $I$. This notation was first introduced for the so-called Stanley-Reisner ideals, which are generated by square-free monomials, see \cite[Chapter 2]{St}.  

\begin{figure}[ht]
\setlength{\unitlength}{0.5cm}
\begin{picture}(16,7)
\put(-5,6){\makebox(0,0){ \bf Fig. 1}}
\put(0,3){\makebox(0,0) {$\begin{array}{ll}I &= (X_1^2X_3^3X_4, X_2^3X_4^2) \\ 
 \sqrt I &= (X_1X_3X_4, X_2X_4)\\ 
\Fcal(\Delta (I)) &= \{ \{1,2,3\}, \{1,4\}, \{3,4\}\} \end{array}$}} 
\hspace{4cm}
\includegraphics[width=7cm,height=4cm]{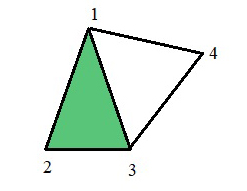} 
 \end{picture}
 \end{figure}
 
 For every $\albf = (\alpha_1,\ldots,\alpha_r) \in \Zset^r$, we define its co-support  to be the set $\Ga := \{i \ | \ \alpha_i < 0\}$.  For a subset $F$ of $[r]$, let $R_F := R[X_i^{-1} \ | \ i \in F]$  be the localization of $R$ at $F$.  
Set
\begin{equation} \label{EQ01}  \Da(I) ;= \{ F \subseteq [r]\setminus \Ga|\  \Xbf^\albf \notin IR_{F\cup \Ga} \}. \end{equation}
We set $\widetilde{H}_i(\emptyset ; K) = 0$ for all $i$,   $\widetilde{H}_i( \{\emptyset \} ; K) = 0$ for all $i\neq -1$, and $\widetilde{H}_{-1} (\{ \emptyset\} ; K) =  K$.   Thanks to \cite[Lemma  1.1]{GH} we may reformulate Takayama's result as follows.

 \begin{thm} \label{Takay} {\rm (\cite[Theorem 2.2]{Ta})} $\dim_K H_{\mfr}^i(R/I)_{\albf} = \dim_K \widetilde{H}_{i-| \Ga|-1}(\Da (I); K).$
\end{thm}

It was  shown in \cite[Lemma 1.3]{MT1} that $\Delta_{\alpha}(I)$ is a subcomplex of $\Delta(I)$.  As a consequence, $ H_{\mfr}^i(R/I)_{\albf} = 0$ provided $\Ga \not\in \Delta (I) $.  Assume that $\Ga \in \Delta (I) $. 
Then, by \cite[Lemma1.2]{HT2},   
$$\Da(I) = \{ F \in \link_{\Delta (I)}(\Ga) |\ x^\alpha \notin IR_{F\cup \Ga} \},$$
 where  the link of  a face $F$ of a simplicial complex $\Delta $ is defined by  
$$\link_\Delta (F) = \{ G \subseteq [r]\setminus F|\ F\cup G \in \Delta \}.$$
 Using this remark, one can see that in the case of Stanley-Reisner ideals, Takayama's formula is exactly the famous Hochster's formula, see \cite[Theorem 4.1]{St}. Hochster's formula plays crucial role in the theory of Stanley-Reisner ideals, where one can find rich interaction between commutative algebra and combinatorics (see, e.g., \cite{BH, St}).  We will see in this paper, that Takayama's theorem is very useful in the study of some invariants of powers of monomial ideals.
 
 Another important notion in commutative algebra is the set of associated primes of a ring.  Now we go back to an arbitrary Noetherian   ring $R$. Recall that a prime ideal  $\pfr \subset  R$ is called an associated prime if $\pfr$ is the annihilator $\ann (x) := \{a\in R|\ ax = 0\}$ for some $x\in R$. The set of associated primes of  $R$ is written as $\ass(R)$. One can say that this set has a central position in commutative algebra like prime divisors of a natural number in number theory. However, it is difficult to compute $\ass(R)$. Therefore the following result by Brodmann is very nice and also finds a lot of application:
 
 \begin{thm} \label{Brod2} {\rm \cite{B1}}  Let $I \subset R$ be a proper ideal. Then the set $\ass(R/I^n)$ is stable for all $n\gg 0$.
 \end{thm}
 
 The research of Brodmann was motivated by a Conjecture of Ratliff in \cite{Ra0} which says that $\ass(R/I^n) \subseteq \ass^\infty (R/I)$ for all $n\ge 1$.  Note that already in \cite[Page 80]{ME}, there is an example constructed to show that the sequence $\{ \ass(R/I^n) \}$ is not monotone.  For integral closures,  extending a result by McAdam and  Eakin \cite[Props. 7 and 18]{ME}, Ratliff shows that a stronger result holds 
 
  \begin{thm} \label{Rat} {\rm \cite[Theorems 2.4 and 2.7]{Ra}} Let $I \subset R$ be a proper ideal. Then the sequence of sets $\{ \ass(R/\overline{I^n}) \}$ is increasing and becomes stable when $n\gg 0$.
 \end{thm}
  
  We would like to know when the sequences $\{ \ass(R/I^n) \}$ and  $\{ \ass(R/\overline{I^n})\}$ become stable. For that, we need
  
  \begin{defn} \label{ast} For a proper ideal $I$ of a Noetherian ring $ R$, set
$$\astab (I) := \min\{ s|\ \ass(R/I^n) = \ass (R/I^s) \ \forall n\ge s \} ,$$
and
$$\astabIt (I) := \min\{ s|\ \ass(R/ \overline{ I^n}) = \ass  (R/ \overline{ I^s}) \ \forall n\ge s \} .$$
 \end{defn}

 \vskip0.3cm
 \noindent {\bf Problem 2}. {\it Give  upper bounds on $\astab(I)$ and $\astabIt(I)$ in terms of other invariants of $R$ and $I$.}
 \vskip0.3cm

In the general case this problem seems to be very hard, because there is no effective way to compute the sets $\ass(R/I^n)$ and  $\ass(R/ \overline{ I^n})$. However, the  prime divisors of a monomial ideal are easily to be found. Therefore one can solve Problem 2 for monomial ideals. This will be summarized in Section \ref{Ass}.

\section{Stability of associated primes} \label{Ass}

From now on, let $R= K[X_1,...,X_r]$, $\mfr = (X_1,...,X_r)$ and $I$ a proper monomial ideal of $R$.   If $r\geq 2$, then for a positive integer $j\leq r$ and $\albf= (\alpha_1,...,\alpha_r) \in \Rset^r$, we set
$$ \albf[j] = (\alpha_1,...,\alpha_{j-1}, \alpha_{j+1},...,\alpha_r).$$
Denote $\Xbf^{\albf}[j]$ the monomial obtained from $\Xbf^{\albf}$ by setting
$X_j=1$. Let $I[j]$ be the ideal of $R$ generated by all monomials $\Xbf^{\albf}[j]$
such that $\Xbf^{\albf}\in I$.  Since any associated prime $\pfr$ of a monomial ideal $\afr$ is generated by a subset of
variables and there is a monomial $m \not\in \afr$ such that $\pfr = \afr: m$, one can easily show

 \begin{lem} \label{AssL1} {\rm \cite[Proposition 4, Lemma 11]{NT1}}  Let $\mfr = (X_1,...,X_r)$ and $r\geq 2$. Then for all
$n\geq 1$ we have:

(i) $\ass(R/I^n) = \ass (I^{n-1}/I^n)$ and $\ass(R/\overline{I^n}) = \ass ( \overline{I^{n-1}}/ \overline{I^n} )$,

(ii) $ \ass(I^n/I^{n+1}) \setminus \{\mfr \} = \cup_{i=1}^r \ass(I[i]^n/I[i]^{n+1}),$

(iii) $ \ass(\overline{I^n}/ \overline{I^{n+1}}) \setminus \{\mfr \} = \cup_{i=1}^r \ass(\overline{I[i]^n}/ \overline{I[i]^{n+1}}).$
\end{lem}
 
\begin{rem} \label{AssRem} On one side, Lemma \ref{AssL1}  allows us to do induction on the number of variables. On the other side, in order to study the stability of the set of associated primes, it reduces to checking if $\mfr \in \ass(R/I^n)$ or $\in \ass(R/\overline{I^n})$, respectively. 
\end{rem}

\subsection{Associated primes of  integral closures of powers}
 
 One can identify a monomial $\Xbf^{\albf}$ with the integer point $\albf \in \Nset^r \subset \Rset^r$. For a subset $A \subseteq R$, the exponent set of $A$ is 
$$E(A) := \{\albf \mid \Xbf^{\albf}\in A\} \subseteq  \Nset^r.$$
 So a monomial ideal $\afr$ is completely defined by its exponent set  $E(\afr)$. Then, we can geometrically describe $\overline{\afr}$ by using its Newton polyhedron.

\begin{defn} \label{NP} Let $\afr$ be a monomial ideal of $R$.  The Newton polyhedron of $\afr$ is $NP(\afr) := \conv\{E(\afr)\}$, the convex hull of the exponent set  $E(\afr)$ of $\afr$ in
the space $\Rset^r$.
\end{defn}

 The following results are well-known (see \cite{RRV}):
$$ E(\overline I) = NP(I) \cap \Nset^r ,$$
and 
\begin{equation}\label{EN2}
NP(I^n) =nNP(I)= n\conv\{E(I)\} +\Rset_{+}^r \ \text{for all}\ n\ge   1.
\end{equation}

The above equalities say that (exponents of)  all monomials  of $\overline{I}$ form the set of  integer points in $NP(I)$ (while we do not know which points among them do not belong to $I$), and the Newton polytope $NP(I^n)$  of $I^n$ is just a multiple of $NP(I)$. 

\vskip0.3cm

\begin{figure}[ht]
\setlength{\unitlength}{0.5cm}
\begin{picture}(16,10)
\put(-5,8){\makebox(0,0){ \bf Fig. 2}}
\put(0,6){\makebox(0,0) {$I = (X_1X_2^4, X_1^3X_2^2, X_1^5X_2)$}} 
\put(0,4){\makebox(0,0){ $\begin{array}{l} \text{The hole means the point} \\
\text{does not belong to}\  E(I). \end{array}$}}
\hspace{4cm}
\includegraphics[width=6cm,height=5cm]{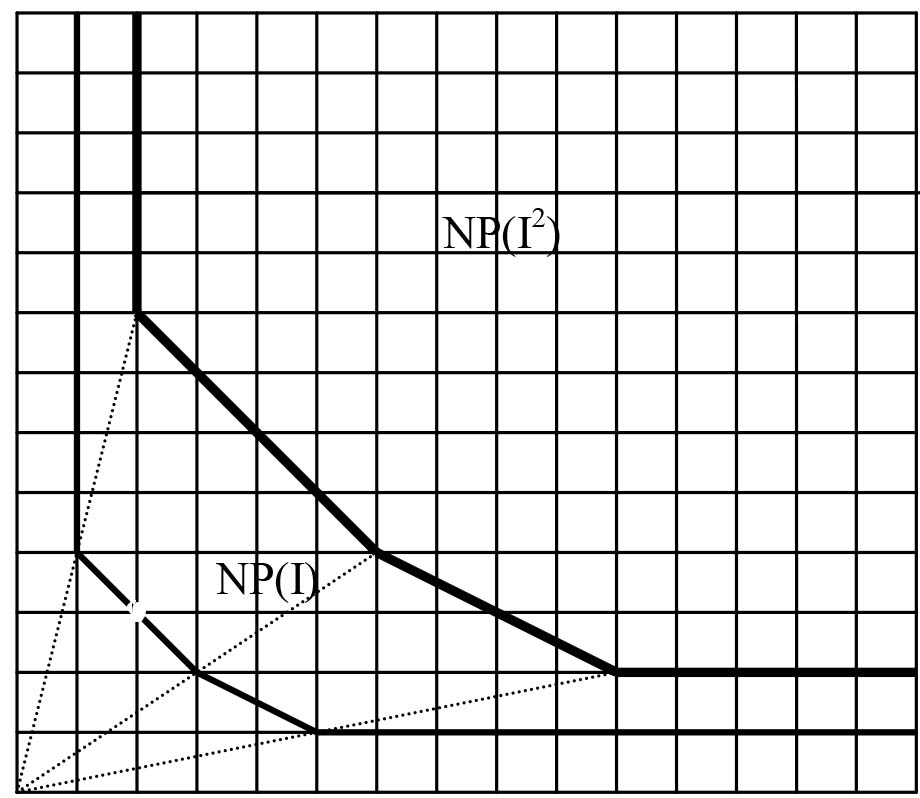} 
 \end{picture}
 \end{figure}

\begin{rem} \label{mAss}
By the definition of $NP(I^n)$ and (\ref{EN2}) it follows that  $\mfr \in \ass(R/\overline{I^n})$ if and only if  there is $\albf \not\in nNP(I)$ and $\albf + \ebf_i \in nNP(I)$ for all $1\le i\le r$, where $\ebf_1,...,\ebf_r$ form the canonical basis of $\Rset^r$.
\end{rem}

Let $G(I)$ denote the minimal monomial generating system of  $I$ and
$$d(I) := \max\{\alpha_1 + \cdots + \alpha_r |\ \Xbf^\albf \in G(I) \},$$
the maximal generating degree  of $I$.    Using convex analysis and lineal algebra, one can show

\begin{lem} {\rm (\cite[Lemma 6]{NT1}, \cite[Lemma 2.2]{HT3})}  \label{NPH} The Newton polyhedron $NP(I)$ is the set of solutions of a system of inequalities of the form
\begin{equation} \label{EN3}
\{\xbf\in\Rset^r \mid \left<\abf_j, \xbf\right> \ge  b_j,\  j=1,\ldots,q\},
\end{equation}
such that each hyperplane with the equation $\left<\abf_j,\xbf\right> = b_j$ defines a facet of $NP(I)$, which contains $s_j$ affinely
independent points of $E(G(I))$ and is parallel to $r - s_j$ vectors of the canonical basis. Furthermore, we can
choose $\zbf \ne \abf_j \in \Nset^r, b_j \in \Nset$ for all $j = 1, . . . , q$; and if we write $\abf_j=(a_{j1},\ldots,a_{jr})$, then
$$a_{ji} \le  s_jd(I)^{s_j-1} \ \text{ for all } i=1,\ldots, r,$$
where $s_j$ is the number of non-zero coordinates of $\abf_j$.
\end{lem}
 
Now one can give an effective necessary condition for $\mfr \in \ass(R/\overline{I^n})$ for some $n>0$. It follows from Remark \ref{mAss} and Lemma \ref{NPH}.

\begin{lem} {\rm \cite[Lemma 13]{NT1}} \label{mAssIn}  Assume that $\mfr \in \ass(R/\overline{I^n})$ for some $n>0$. Then there is a vector $\abf_i$ determined  in Definition \ref{EN3} such that $\abf_i > 0$, that is $a_{ij} > 0$ for all $1\le j\le r$.
\end{lem}

\begin{lem} \label{AssL2} Let $I$ be a monomial ideal in $R$ with $r > 2$. If $\mfr \in \ass (R/\overline{I^s})$ for some $s\ge  1$, then
$\mfr \in \ass (R/\overline{I^n})$ for all $ n\ge  (r-1)r d(I)^{r-2}.$
\end{lem}

\begin{proof} (Sketch):  Let $m := (r-1)rd(I)^{r-2}$. Since the sequence $\{\ass (R/\overline{I^n})\}_{n\ge  1}$ is increasing by  Theorem 
\ref{Rat},  it suffices to show that $\mfr \in \ass (R/\overline{I^m})$.  As $\mfr \in \ass (R/\overline{I^s})$, by Lemma  \ref{mAssIn}, there is a supporting hyperplane of $NP(I)$, say $H$,  of the form $\left<\abf, \xbf \right> = b$ such that all coordinates of $\abf$ are positive. By Lemma  \ref{NPH}, this hyperplane passes through $r$ affinely independent points of $E(G(I))$, say $\albf_1, \ldots ,\albf_r$. Denote the barycenter  of the simplex $[\albf_1, \ldots ,\albf_r]$ by $\albf := \frac{1}{r} (\albf_1 + \cdots + \albf_r)$, and  let $\bebf:= m\albf -\ebf_r$. Then one can show that $\bebf \not\in mNP(I)$ and $\bebf + \ebf_i \in  mNP(I)$ for all $i\ge 1$. Hence  the lemma follows from  Remark \ref{mAss}.
\end{proof} 

The first bound on $\astabIt$  is given in \cite[Theorem 16]{NT1}. It is then improved as follows.

\begin{thm} {\rm \cite[Theorem 2.7]{HT3}}  \label{AssTh1} Let $I$ be a monomial ideal of $R$. Then
\begin{equation*} \astabIt(I) \le 
\begin{cases}
1 & \text{ if } \ell(I) \le  2,\\
\ell(I)(\ell(I)-1)d(I)^{\ell(I)-2} & \text{ if } \ell(I) > 2.
\end{cases}
\end{equation*}
\end{thm}

This theorem almost immediately follows from Lemma \ref{mAssIn} and Remark  \ref{AssRem} by using induction on $r$ (based on Lemma 
 \ref{AssL1}).

\begin{rem} By \cite[Theorem 2.3]{BA}, we can compute $\ell(I)$ in terms of geometry of  $NP(I)$.
$$\ell(I) = \max\{\dim F + 1 \mid F \text{ is a compact face of } NP(I)\}.$$
\end{rem}

\begin{exm} \cite[Prposition 17]{NT1} \label{Example1}
  Let $r\geq 4$ and $d>r-3$.
We put
$$u = X_1^{r-3 \choose 0} X_2^{r-3 \choose 1}\cdots  X_{r-3}^{r-4 \choose 0} \ \text{and} \
v= X_1^{\beta_1}\cdots  X_{r-3}^{\beta_{r-4}}X_{r-2}^{d-r+2},$$
where
$$\beta_i =   \begin{cases}  0 & \text{if}\ r-3-i\  \text{is\ even},\\
2{r-3\choose i} & \text{if}\ r-3-i\  \text{is\ odd}.
 \end{cases}$$
Let
$$I= (uX_1^d,uX_2^{d-1}X_r,...,uX_{r-2}^{d-r+3}X_r^{r-3}, uX_{r-1}X_r^{d-1}, vX_r^{r-3}).$$
It is generated by monomials of the same degree  $d(I) = d+2^{r-3} -1$.   Then  $\mfr \in \ass(R/\overline{I^n})$ for all $n\gg 0$ and if $\mfr \in \ass(R/\overline{I^n})$, then $n\ge n_0 := \frac{d(d-1)\cdots (d-r+3)}{r(r-3)}$. In particular,
$$\astabIt(I) \ge  n_0.$$
This shows that the bound in Theorem  \ref{AssTh1}  is almost optimal, and {\it it must  depend on the maximal generating degree of $I$}. \end{exm}

\begin{proof}(Sketch):  In this example, by Lemma \ref{mAssIn}, $\mfr\in \ass(R/\overline{I^n})$ for $n\gg 0$. On the other hand,  the projection of $\conv(E(G(I))$ into the hyperplane $\Rset^{r-1}$ of the first $(r-1)$ coordinates form a simplex, say $\Delta $. Using Lemma 
 \ref{NPH} and Remark \ref{mAss}, one can show that if $\mfr \in  \ass(R/\overline{I^n})$ for a fix $n$, then the simplex $n\Delta $ must contain $r$ integer points of the form $\bebf',\ \bebf' + \ebf'_1, ...., \bebf' + \ebf'_{r-1} \subset \Nset^{r-1}$, where $\ebf'_1,...,\ebf'_{r-1}$ are unit vectors of $\Rset^{r-1}$. The simplex $\Delta $ is so far from  being  ``regular", that only its big multiples satisfy this combinatorial property.
\end{proof}

\vskip0.3cm
\noindent {\bf Question 3}. Assume that $I$ is a square-free monomial ideal. Is there a linear upper bound on $\astabIt(I)$ in term of $r$?

\subsection{Associated primes of powers}

In the sequel, by abuse of terminology, for a linear functional
$$\varphi (\xbf) = a_1x_1+ \cdots + a_rx_r,$$
where $a_i\in \Rset$,  we say that $\varphi(\xbf) \geq 0$ is a homogeneous linear constraint, while $\varphi(\xbf) \geq b$ is a linear
constraint. Unlike  integer closures, it is  much more difficult to describe the set of monomials in $I^n$ by linear constrains. However, we have the following observation: 

Assume that the  monomials $\Xbf^{\albf_1},..., \Xbf^{\albf_s}$ generate the  ideal $I$. Then a monomial $\Xbf^{\albf}\in I^m$ if
and only if there are nonnegative integers $a_1,...,a_{s-1}$,  such that $m \geq a_1+\cdots + a_{s-1}$ and $\Xbf^{\albf}$ is divisible by
$$(\Xbf^{\albf_1})^{a_1}\cdots (\Xbf^{\albf_{s-1}})^{a_{s-1}}
(\Xbf^{\albf_s})^{m- a_1- \cdots - a_{s-1}}.$$
This is equivalent to
$$\alpha_j \geq \alpha_{1j}a_1 + \cdots + \alpha_{(s-1)j}a_{s-1} + \alpha_{sj}(m- a_1-\cdots -
a_{s-1}),$$ for all $j=1,...,r$.

From this observation, $\albf \in E(I^m)$ if and only is it is a part of an integer solution of a system of linear constrains in $r+s-1$ variables. Unfortunately this correspondence is not one-to-one, so that we cannot reverse a constrain in order to get a criterion for $\albf \not\in E(I^m)$. Nevertheless this observation is  useful in finding an upper bound on $\astab(I)$ in \cite{Hoa}. 

The next observation is that in this case thanks to Lemma \ref{AssL1}(i), it is easier to work with $\ass (I^{n-1}/I^n)$ than with $\ass(R/I^n)$, because the quotient modules $I^{n-1}/I^n$, $n\ge 1$, can be put together in the so-called associated graded ring of $I$:
$$G = \oplus_{n\geq 0}I^n/I^{n+1}.$$
 Further, $ \mfr \in \ass(I^n/I^{n+1})$
if and only if the local cohomology module $H^0_{\mfr}(I^n/I^{n+1}) \neq 0$. This local cohomology can be computed as follows (see  \cite[Lemma 3.2]{Hoa})
$$H^0_{\mfr G}(G)_{n-1} \cong H^0_{\mfr}(I^{n-1}/I^n) \cong
 \frac{I^{n-1} \cap I[1]^n \cap \cdots \cap I[r]^n}{I^n}.$$
 
 Using the above observation, one can associate the family of $E(I^{n-1} \cap I[1]^n \cap \cdots \cap I[r]^n$) to a set $\Ecal  \subset \Nset^{rs +s}$ of  integer solutions of linear constrains in $rs +s$ variables. If we denote the set of integer solutions of the corresponding system of homogeneous linear constrains by $\Scal$, then $\Scal$ is a semigroup, so that $K[\Scal]$ is a ring, and $K[\Ecal]$ is a $K[\Scal]$-module. One can prove that $H^0_{\mfr G}(G)$ is isomorphic to a quotient of $K[\Ecal]$ (\cite[Lemma 3.4]{Hoa}). Using linear algebra and Caratheodory's Theorem (see, e.g. \cite[Corollary 7.1(i)]{Sch}) one can show that the maximal generating degree of $K[\Ecal]$ over $K[\Scal]$ is bounded by 
 $$B_1 := d(rs+s+d)(\sqrt{r})^{r+1}(\sqrt{2} d)^{(r+1)(s-1)},$$
 where $d= d(I)$, $s$ the number of minimal generators of $I$ (see \cite[Proposition 3.1]{Hoa}).
From that one obtains

\begin{prop} {\rm \cite[Proposition 3.2]{Hoa}} \label{D5}  Let $n\geq B_1$ be an integer. Then
$$\ass(I^n /I^{n+1}) \supseteq \ass(I^{n+1}/I^{n+2}).$$
\end{prop}

In order to get the reverse inclusion we use another local cohomology module. Recall that the  Rees algebra $\Rcal := \Rcal(I) = \oplus_{n\ge 0}I^nt^n$. Let $\Rcal_+ = \oplus_{n>0}I^nt^n$.  
The local cohomology module $H^0_{\Rcal_+}(G)$ is also a $\Zset$-graded
$\Rcal$-module. Let
$$a_0(G) = \sup\{n | \ H^0_{\Rcal_+}(G)_n \neq 0 \}.$$
(This number is to be taken as $-\infty$ if $H^0_{\Rcal_+}(G)=0$.) It is
related to the so-called  Castelnuovo-Mumford regularity of $G$
(see, e.g., \cite{Sh}). Then   S. McAdam and P.
Eakin  show that
 $ \ass(I^n /I^{n+1}) \subseteq
\ass(I^{n+1}/I^{n+2})$ for all $n> a_0(G)$ (see   \cite[pp. 71, 72]{ME},  and also \cite[Proposition 2.4]{Sh}).
Now one can again use (another) system of linear constrains  to show  

\begin{prop} {\rm \cite[Proposition 3.3]{Hoa}} \label{D9} We have
$$a_0(G) < B_2 := s(s+r)^4 s^{r+2} d^2(2d^2)^{s^2-s+1}.$$
\end{prop}

Putting together Propositions \ref{D5} and \ref{D9}, we get 

\begin{thm} {\rm \cite[Theorem 3.1]{Hoa}}  \label{AssTh2} We have
$$\astab(I) \le \max\{d(rs+s+d)(\sqrt{r})^{r+1}(\sqrt{2} d)^{(r+1)(s-1)}, \
s(s+r)^4 s^{r+2} d^2(2d^2)^{s^2-s+1} \}.$$
\end{thm}

Of course, one can bound $s$ in terms of $d$ and $r$. But then  the resulted bound would be a double exponential bound.
In spite of Theorem \ref{AssTh1}, we would like to ask:

\vskip0.3cm
\noindent {\bf Question 4}. {\it i) Is there an upper bound on $\astab(I)$ of the order $d(I)^r$?

ii) Assume that $I$ is a square-free monomial ideal. Is there a linear upper bound on $\astab(I)$ in term of $r$? } 
\vskip0.3cm

No that using Example 1, one can construct an example to show that in the worst case, an upper bound on $\astab(I)$ much be at least of the order $d(I)^{r-2}$ (provided that $r$ is fixed), see \cite[Example 3.1]{Hoa}. 

Another interesting problem is to find the stable set $\ass^\infty (R/I)$ for $n\gg 0$. There is not much progress in this direction. However,  Bayati, Herzog and Rinaldo  can completely solve a kind of reverse problem: In \cite{BHR} they prove that any set of nonzero monomial prime ideals can be realized as the stable set of associated primes of a monomial ideal.

\subsection{Case of edge ideals and square-free monomial ideals}

There are some partial classes of monomial ideals where $\astab(I)$ is bounded by a linear function of $r$, that gives a partial affirmative answer to Question 4. The best results are obtained for square-free monomial ideals whose generators are all of degree two. Then one can associate such an ideal to a  graph. More precisely, let $G = (V,E)$ be a simple undirected graph with the vertex set $V = [r] := \{1,2,...,r\}$ and the edge set $E$. The ideal
$$I(G) = (x_ix_j|\ \{i,j\} \in E) \subset R,$$
is called  {\it edge ideal}.  Recall that $G$ is {\it bipartite} if $V= V_1\cup V_2$ such that $V_1\cap V_2 = \emptyset $ and  there is no edge connecting two vertices of the same set $V_i$. In our terminology, one can reformulate \cite[Theorem 5.9]{SVV} as follows.

\begin{thm} \label{SVVThm} The graph $G$ is bipartite if and only if $\astab (I(G)) = 1$.
\end{thm}

The following result is a reformulation of \cite[Corollary 2.2]{CMS}, which reduces the problem of bounding $\astab(I(G))$ to the case of connected graphs. 

\begin{lem} \label{CMSlem} Assume that $G_1,...,G_s$ are connected components of $G$. Then
$$\astab(I(G)) \le \sum_{i=1}^s \astab(I(G_i)) - s + 1.$$
\end{lem}

Recall that a cycle $C$ in a graph is a sequence of different vertices $\{ i_1,...,i_s\} \subset V$ such that $\{i_j, i_{j+1}\} \in E$ for all $j\le s$, where $i_{s+1} \equiv i_1$. The number $s$ is called the length of  $C$. It is an elementary  fact in the graph theory that $G$ is bipartite if and only if $G$ does not contain  odd cycles.

\begin{figure}[ht]
\setlength{\unitlength}{0.5cm}
\begin{picture}(14,12)
\put(-5,6){\makebox(0,0){ \bf Fig. 3}}
\includegraphics[width=10cm,height=6cm]{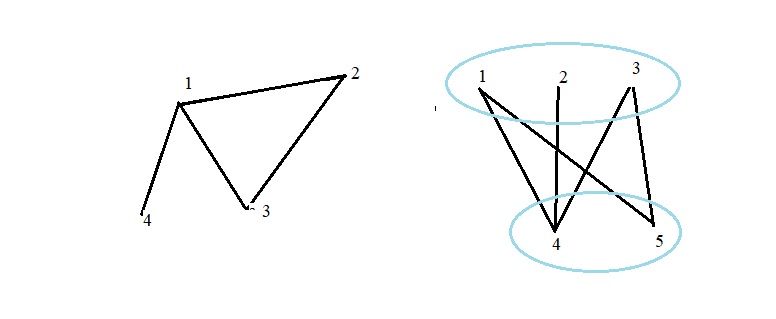}
\put(-10,0.5){\makebox(0,0){Non-bipartite graph $G_1$  \hskip1.5cm    Bipartite graph $G_2$}}  
 \end{picture}
 \end{figure}
 
 The following result together with Theorem \ref{SVVThm} and Lemma \ref{CMSlem} gives a good bound on $\astab(I(G))$.
 
 \begin{thm} {\rm \cite[Lemma 3.1, Proposition 4.2]{CMS}} \label{CMSThm} Assume that  $G$ is non-bipartite connected graph. Let $2k+1$ be the smallest odd cycle contained in $G$ ($k\ge 1$). 
 
 (i) If $G$ is a cycle (of length $2k+1$), then $\astab(I(G)) = k+1$,
 
 (ii) If $G$ is not a cycle, then $\astab(I(G)) \le  r-k-1$.
\end{thm}

The proofs of the above results in \cite{CMS} are quite long, but they are of combinatorial nature and do not require much from the graph theory. Even the above bound is quite good, by using other invariants, one can get a better bound. For an example, in the situation of Theorem \ref{CMSThm}(ii), if $G$ has $s$ vertices of degree one, then  $\astab(I(G)) \le  r-k-s$ (see \cite[Corollary 4.3]{CMS}).

If we set $\Min(I)$ to be the set of minimal associated primes of $R/I$, then $\Min(I) \subseteq \ass(R/I^n)$ for all $n\ge 1$. An element from $\ass(R/I^n)\setminus \Min(I)$ is called an {\it embedded associated prime}. Therefore, in order to study the stability of $\ass(R/I^n)$, it suffices to study embedded associated primes.

In \cite{HiLT}, using the Takayama's Theorem \ref{Takay}, Hien, Lam and N. V. Trung  show that embedded primes of $R/I(G)^n$ are characterized by the existence of (vertex) weighted graphs with special matching properties, see  \cite[Theorem 2.4]{HiLT}.  There are (infinite) many weighted graphs with the base graph $G$. Using techniques in the graph theory, they  can give  necessary or sufficient conditions for an ideal generated by a subset of variables to be an embedded associated primes of $R/I(G)^n$ in terms of vertex covers of $G$ which contain certain types of subgraphs of $G$, see  \cite[Theorem  2.10, Theorem 3.5]{HiLT}. From that they  derive a stronger bound on $\astab(I(G))$, see \cite[Corollary 3.7]{HiLT} and also \cite[Example 3.8]{HiLT}.  Moreover, their method gives an algorithm to compute $\ass(R/I(G)^t)$ for a fixed integer $t$. This was done for $t=2,3$ in \cite{HiLT} and for $t=4$ in \cite{HiLam}. 

In \cite{LamT},  after extending  some results in the graph theory,  Lam and N. V. Trung are able to characterize  the existence of  weighted graphs with special matching properties in terms of the so-called generalized ear decompositions of the base graph $G$. 
Using  this notion they  give a new upper bound on $\astab(I(G))$ (see \cite[Theorem 4.7]{LamT}). Moreover, they can give a precise formula for $\astab(I(G))$.   In order to formulate  their results one needs to introduce some rather technical notions on graphs. Therefore we do not go to the details here. We only want to state the following nice main result of \cite{MMoV},  that they  can derive as an immediate consequence of their results.

\begin{thm} {\rm (\cite[Theorem 2.15]{MMoV}, \cite[Corollary 4.5]{LamT})}. \label{AssAscen} The sequence $\{\ass (R/I(G)^n)\}$ is ascending, that is $\ass(R/I(G)^n) \subseteq \ass(R/I(G)^{n+1})$ for all $n\ge 1$.
\end{thm}

This property does not hold for square-free monomial ideals. Using a counter-example in combinatorics, Kaiser et al. constructed a square-free monomial ideal $J$ in 12 variables such that $\ass(R/J^4)$ is a proper subset of $\ass(R/J^3)$, see \cite[Theorem 11]{KSS}. 

For an arbitrary monomial $I$, the sequence $\{\ass (R/I^n)\}$ is even not necessarily monotone. The first example is given in  page 549 of \cite{HH2}. Recently, there is given a general construction of monomials ideals $I$ for which the non-monotonicity of $\{\ass (R/I^n)\}$ can be arbitrarily long.

\begin{thm} {\rm \cite[Corollary 6.8]{HHTT}}  \label{AssNotMon} Let $A$ be any finite set of positive integers. Then there exists a monomial ideal $I$ in a polynomial ring $R$ such that $\mfr \in \ass(R/I^n)$  if and only if $n \in A$. 
\end{thm}
  
The ascending property of Ass, that means $\ass(R/I^n) \subseteq \ass(R/I^{n+1}), \ n\ge 1$, is also referred as the {\it persistence property} (with respect to associated ideals).  This property is important in finding the stable set $\ass^\infty (R/I)$.

For arbitrary square-free monomial ideals there are some interesting results given in \cite{FHaV1, FHaV2, HaMo, HQ}, that  relate $\astab(I)$   to combinatorics.  In order to say about their study we need to introduce some notions. 

The set of bases of a polymatroid of rank $d$ based on $[r]$ is a set $\Bcal$ of integer points $\albf \in \Nset^r$ satisfying the following conditions:

\begin{itemize}
\item $|\albf| := \alpha_1 + \cdots + \alpha_r = d$ for all $\albf \in \Bcal$,
\item   (Exchange property) For all $\albf, \bebf \in \Bcal$ for which $\alpha_i > \beta_i$ for some {i}, there exists $j$ such that $\beta_j > \alpha_j$ and $\albf - \ebf_i + \ebf_j \in \Bcal$.
\end{itemize}

A monomial ideal $I$ is called a {\it polymatroidal ideal}, if there exists a set of bases $\Bcal \subset \Zset^r$ of a polymatroid, such that 
$I = (\Xbf^\albf|\ \albf \in \Bcal)$.

\begin{thm} {\rm  \cite[Proposition 2.4 and Theorem 4.1]{HQ}} Let $I$ be a polymatroidal ideal. Then $I$ satisfies persistence property and $\astab(I) \le \ell(I)$.
\end{thm}

A  simple hypergraph $\Hcal$ is a pair of  the vertex set $V = [r]$  and an edge set $\Ecal = \{E_1,...,E_t\}$, where $E_i \subseteq V$. We assume that $\Hcal$ has no isolated vertices, each $E_i$ has at least two elements and that  $E_i \not\subseteq E_j$ for all $i\neq j$. When the $E_i's$ all have cardinality two, then $\Hcal$ is a simple graph.  The ideal generated by all square-free monomials $X_{i_1}\cdots X_{i_s}, \ \{i_1,...,i_s\} = E_i$, is called {\it edge ideal} of $\Hcal$ and denoted by $I(\Hcal)$. For a very special class of $\Hcal$,  Ha and Morey determine the least number $k$ such that $\mfr \in \ass(R/I^k)$ and $\mfr \not\in \ass(R/I^t)$ for all $t<k$, see \cite[Theorem 4.6]{HaMo}. 

A vertex cover of $\Hcal$  is a subset $W$ of  $V$ such that if $E\in \Ecal$, then $W\cap E \neq \emptyset $. A vertex cover is minimal if no proper subset is also a vertex cover. Denote $J = J(\Hcal)$  the {\it cover ideal} of $\Hcal$, which is generated by the square-free monomials corresponding to the minimal vertex covers of $\Hcal$. Francisco, Ha and Van Tuyl propose a conjecture related to the chromatic number of a graph $G$, and prove the persistence property of $\ass(R/J(G)^n)$ provided that the conjecture holds (see \cite[Theorem 2.6]{FHaV1}). In another paper, they give  an explicit description of all associated primes of $\ass(R/J(\Hcal)^n)$, for any fixed number $n\ge 1$, in terms of the coloring properties of hypergraphs arising from $\Hcal$, see \cite[Corollary 4.5]{FHaV2}. From this description they give a lower bound on $\astab(J(\Hcal))$.

Recall that a $t$-coloring of $\Hcal$  is any partition of $V = C_1\cup \cdots \cup C_t$  into $t$ disjoint sets such that for every $E\in \Ecal$, we have $E \not\subseteq C_i$ for all $i= 1,...,d$.  The $C_i's$ are called the {\it color classes}.  The chromatic number of $\Hcal$, denoted $\chi(\Hcal)$, is the minimal $t$ such that $\Hcal$ has a $t$-coloring.

\begin{prop} {\rm \cite[Corollary 4.9]{FHaV2}} $\astab(J(\Hcal)) \ge \chi (\Hcal) -1$.
\end{prop}

Inspired by this result they pose the following question:

\vskip0.3cm
\noindent {\bf Question 5} \cite[Question 4.10]{FHaV2}. {\it For each integer $\ge 0$, does there exist a hypergraph $\Hcal_n$ such that the stabilization of associated primes occurs at $a \ge \chi (\Hcal_n) - 1) + n$? }
\vskip0.3cm

As one can see from the above discussion, all results concerning the existence of a linear bounding on $\astab(I)$ and the persistence property are given for very special squaree-free monomial ideals. Nevertheless, these results  establish  surprising relationships between some seemingly unrelated notions of commutative algebra and combinatorics and raise many more problems and questions.  Thus, they will stimulate intensive investigation in the near future.

\section{Stability of Depth} \label{Depth}

\subsection{Depth of powers of integral closures}

Due to some reasons, we can completely solve Problem 1 for $\dstabIt(I)$. First,  as an immediate consequence of Theorem \ref{Takay}, one can get  the  following ``quasi-decreasing" property of the depth function $\depth  R/\overline{I^n}$. We don't know if this property holds for integral closures of powers of an arbitrary homogeneous  ideal.

\begin{lem} {\rm \cite[Lemma 2.5]{HT3}} \label{DepthInc} For any monomial ideal $I$ of $R$, we have
\begin{enumerate}
\item $\depth R/\overline{I^m}\ge  \depth R/\overline{I^{mn}} \text{ for all } m, n\ge 1$.
\medskip
\item $\lim_{n\rightarrow \infty}\depth R/\overline{I^n} =\dim R-\ell(I),$ where $\ell (I)$ denotes the analytic spread of $I$.
\end{enumerate}
\end{lem}

\begin{proof} (Sketch): For proving the first statement one can set $m=1$.  Then, using Theorem \ref{Takay} and the fact that for any $\albf \in \Zset^r$, $CS_{n\albf}= \Ga$ and $\Da(\overline{I}) = \Delta_{n\albf}(\overline {I^n})$, one can quickly show that $H^t_{\mfr}(R/\overline{I})_{\albf} \ne 0$ implies $H^t_{\mfr}(R/\overline{I^n})_{n\albf} \ne 0$.

The second statement follows from \cite[Propostion 3.3]{EH} and the fact that $\overline{I^{r-1}}$  is torsion-free (see \cite[Corollar 7.60]{Vs}.
\end{proof} 

\vskip0.3cm
\noindent {\bf Question 6}. {\it  Is the depth function $\depth R/\overline{I^m}$ decreasing?}
\vskip0.3cm

As we can see from Theorem  \ref{Takay}, in order to study the local cohomology module, we need to have an effective description of  $\Da(\overline{I^n})$. In the case of integral closures, we do have it. Keeping the notations in Lemma \ref{NPH}, we set $\supp(\abf_j) := \{i \mid a_{ji} \ne 0\}$.  

\begin{lem} {\rm \cite[Lemma 3.1]{HT3}} \label{FNPn} For any $\albf \in \Nset^r$ and $n\ge  1$, we have
$$\Da(\overline{I^n}) =\left< [r]\setminus \supp(\abf_j) \mid j\in \{1,\ldots,q\} \text{ and }  \left<\abf_j,\albf \right> <nb_j  \right>.$$
\end{lem}

The following lemma is the main step in the proof of  Theorem \ref{DStab}.
 
 \begin{lem} {\rm \cite[Lemma 3.2]{HT3}} \label{DStabL} Let $m\ge  1$ and $t:=\depth R/\overline{I^m}$. Assume that $H_{\mfr}^t(R/\overline{I^m})_{\bebf} \ne 0$ for some $\bebf\in \Nset^r$. If $r\ge  3$, then 
$$\depth R/\overline{I^n} \le t \ \text{ for all } n\ge  r(r^2 - 1)r^{r/2}(r-1)^rd(I)^{(r-2)(r+1)}.$$
\end{lem}
\begin{proof} (Sketch): Assume that
$$\left<\abf_j, \bebf \right> < b_j \text{ for } j =1,\ldots, p, $$
and
$$\left<\abf_j, \bebf \right> \ge  b_j  \text{ for } j  =p+1,\ldots, q,$$
for some $0\le  p\le  q$. Then, by Lemma $\ref{FNPn}$,
$$\Delta_{\bebf}(\overline{I^m}) = \left<[r]\setminus \supp(\abf_j) \mid j=1,\ldots,p\right>  .$$
For each $n \ge  1$, put
$$\Gamma(\overline{I^n}) := \{\albf \in \Nset^r \mid \Delta_{\albf}(\overline{I^n}) = \Delta_{\bebf}(\overline{I^m})\},$$
and 
\begin{equation}\label{EDStabL2}
C_n:=\{\xbf \in\Rset_+^r\mid \left<\abf_j, \xbf \right> < nb_j,  \left<\abf_l, \xbf \right> \ge   nb_l \text{ for } j \le  p; \  p+1\le l \le  q\} \subseteq \Rset_+^r.
\end{equation}
Assume that  $C_n \cap \Nset^r  \neq \emptyset$.  Then  for any $\albf\in C_n \cap \Nset^r$,  by Theorem  \ref{Takay}, we will have
$ H_{\mfr}^t(R/\overline{I^n}))_{\albf} \ne 0,$ whence  $\depth R/ \overline{I^n} \le  t$. It remains to show that $C_n \cap \Nset^r \ne \emptyset$ for any $n\ge r(r^2 - 1)r^{r/2}(r-1)^rd(I)^{(r-2)(r+1)}$.
\end{proof}

\begin{rem} \label{Ehr} From the above sketch of proof we can see that the main technique is to find a number $n_0$ such that $C_n$ contains an integer point for all $n\ge n_0$, or equivalently  that the system of linear constrains in (\ref{EDStabL2}) do have integer solutions. This is  related to the research carried  out by Ehrhart in \cite{E1, E2}, where he shows that the number of integer points in the closure $\overline{C_n} \subset \Rset^r$ is a quasi-polynomial!
\end{rem}

Using Lemma \ref{DStabL}, Theorem  \ref{Takay} and induction on $r$, one can prove the following main result of this subsection.

\begin{thm} {\rm  \cite[Theorem 3.3]{HT3}} \label{DStab}  Let $I$ be a monomial ideal of $R$.  Then
 \begin{equation*}\dstabIt(I) \le 
\begin{cases}
1 & \text{ if } r \le  2,\\
r(r^2 - 1)r^{r/2}(r-1)^r d(I)^{(r-2)(r+1)}  & \text{ if } r > 2.
\end{cases}
\end{equation*}
\end{thm}

 It seems that this bound is too big. However, Example \ref{Example1} shows that an upper bound on $ \overline{\dstab}(I)$ must be at least of the order $d(I)^{r-2}$. 

\vskip0.3cm
\noindent {\bf Question 7}. {\it Is $\dstabIt(I)$ is bounded by a function of the order $d(I)^r$?}
\vskip0.3cm

\subsection{Depth of symbolic powers }

The $n$-th symbolic power of  an ideal $\afr \subset R= K[X_1,...,X_r]$ is the ideal
$$\afr^{(n)} := R\cap(\cap_{\pfr\in \Min(I)} \afr^nR_{\pfr}).$$ 
In other words, $\afr^{(n)}$ is the intersection of the primary components of $\afr^n$ associated to the minimal primes of $\afr$. 

 When $K$ is  algebraically closed  and $\afr$ is a radical ideal,  Nagata and Zariski showed that $\afr^{(n)}$ consists of polynomials in $R$ whose partial derivatives of orders up to $n-1$ vanish on the zero set of $\afr$. Therefore, symbolic powers of an ideal carry richer geometric structures and more subtle information than ordinary powers!
 
 Unlike the ordinary powers, the behavior of $\depth(R/I^{(n)})$ is much more mysterious. If $I$ is a monomial ideal, then the symbolic Rees algebra $\mathcal R_s(I) = \oplus_{n\ge 0}I^{(n)}$ is finitely generated (see \cite [Theorem 3.2]{HHT}). Then Brodmann's Theorem \ref{Brod1}  implies that $\depth(R/I^{(n)})$ is periodically constant for $n\gg 0$.  Very recently D. H. Nguyen and N. V. Trung are able to construct  monomial ideals $I$ for which $\depth(R/I^{(n)})$ is not constant for $n\gg 0$, see \cite[Theorem 5.4]{HopT}. Moreover, they can show
 
\begin{thm} {\rm \cite[Theorem 5.1]{HopT}} \label{HopTThm} Let $\varphi (t)$  be an arbitrary asymptotically periodic positive numerical
function. Given a field $K$, there exist a polynomial ring $R$ over a purely transcendental  extension of $K$ and a homogeneous ideal $I \subset R$ such that $\depth(R/I^{(t)}) = \varphi (t)$ for all $t\ge 1$.
\end{thm}

The construction in \cite{HopT} only gives non monomial ideals. Therefore we would like to ask:

 \vskip0.3cm
\noindent {\bf  Question 8}. {\it   Does Theorem  \ref{HopTThm} hold for the class of monomial ideals?}
\vskip0.3cm

From now until the rest of this subsection, assume that $I$ is a square-free monomial ideal. Let $\Delta = \Delta (I)$. The correspondence $I \leftrightarrow \Delta (I)$ is one-to-one, and we also write $I_\Delta $ for $I$. Assume that $\mathcal F(\Delta) = \{F_1, \ldots, F_m\}$. Then 
$$I_{\Delta} = \bigcap_{F\in \mathcal F(\Delta)} P_F,$$
where $P_F$ is the prime ideal of $R$ generated by variables $X_i$ with  $i\notin F$, and
$$I_{\Delta}^{(n)} = \bigcap_{F\in \mathcal F(\Delta)} P_F^n.$$
In this case, it immediately follows from \cite[Theorem 4.7]{HT1} that $\depth (R/I^{(n)})$ is constant for all $n\gg 0$. Hence, for a Stanley-Reisner ideal $I_\Delta $, we can introduce the following notation.
 
\begin{defn} Let $I$ be a square-free monomial ideal $I$. Set
$$\dstab^*(I) := \min\{m \ge  1 \mid \depth (R/I^{(n)}) = \depth (R/I^{(m)}) \text{\ for all\ } n \ge  m\}.$$
 \end{defn}
 We can define the {\it symbolic analytic spread} of $I_{\Delta}$ by 
$$\ell_s(I_{\Delta}) := \dim \mathcal R_s(I_{\Delta})/\mfr \mathcal R_s(I_{\Delta}).$$
 Let $\bight(I_{\Delta})$ be the {\it big height} of $I_{\Delta}$. 

\begin{thm} {\rm \cite[Theorem 2.4]{HKTT}} \label{SPThm1}  Let $I_{\Delta}$ be a Stanley-Reisner ideal of  $R = k[X_1, \ldots , X_r]$. Then:
\begin{enumerate}
\item $\depth ( R/I_{\Delta}^{(n)}) \ge \dim R -\ell_s(I_{\Delta})$ for all $n \ge  1$;
\item $\depth (R/I_{\Delta}^{(n)}) = \dim R - \ell_s(I_{\Delta})$ for all $n \ge  r(r+1)\bight(I_{\Delta})^{r/2}$.
\end{enumerate}
In particular, $\dstab^*(I_\Delta ) \le r(r+1)\bight(I_{\Delta})^{r/2}$.
\end{thm}

The idea of the proof of this theorem is similar to that of Theorem \ref{DStab}, because in this case we can also effectively compute $\Delta_{\albf}(I_{\Delta}^{(n)})$.

\begin{lem} {\rm (\cite[Lemma 1.3]{MT1})} \label{L02} For all $\albf\in \Nset^r$ and $n \ge 1$, we have
$$\Delta_{\albf}(I_{\Delta}^{(n)}) = \left < F\in \mathcal F(\Delta) \ | \ \sum_{i\notin F} \alpha_i \le n-1\right>.$$
\end{lem}

We think that this bound is too big. Therefore, we would like to ask

\vskip0.3cm
\noindent {\bf Question 9}. {\it Assume that $I$ is a square-free monomial  ideal $I$.

(i) Is the depth function $\depth (R/I^{(n)})$ decreasing?

(ii) Is there a linear bound on $\dstab^*(I) $ in terms of $r$?}
\vskip0.3cm

Note that the ``quasi-decreasing property" of  $\depth (R/I^{(n)})$ can be proved similarly to  Lemma  \ref{DepthInc}. There are some partial positive answers to this question. As a corollary of \cite[Theorems 2.3 and 2.4]{MT1} (also see \cite[Lemma 2.1]{HT2}), we get

\begin{prop} Assume that $\dim R/I \le 2$. Then $\depth (R/I^{(n)})$ is decreasing and $\dstab^*(I) \in \{1,2,3\}$.
\end{prop}

Let $G$ be a simple graph with the vertex set $V = [r]$. Then the cover ideal 
 $$J(G) = \bigcap_{\{i,j\}\in E(G)} (X_i,X_j) \subset R .$$
  It is clear that every unmixed squarefree  monomial ideal of height two is uniquely correspondent to a cover ideal of a graph and vice verse.  Using the so-called polarization technique one can show the non-increasing property of $\depth ( R/J(G)^{(n)})$. Note that this property does not hold  for the sequence $\depth (R/J(G)^n)$ on ordinary powers of $J(G)$ (see \cite[Theorem 13]{KSS}). The graph constructed there has 12 vertices and $\depth(R/J(G)^3) = 0$ while $\depth(R/J(G)^4) = 4$.

\begin{thm} {\rm \cite[Theorem 3.2]{HKTT}}\label{Decreasing} Let $G$ be a simple graph. Then for $n \ge 2$,
$$\depth (R/J(G)^{(n)}) \le \depth (R/J(G)^{(n-1)}) .$$
\end{thm}

 In order to formulate an effective bound on  $\dstab^*(J(G))$ we recall some terminology from the graph theory. A set $M \subseteq E(G)$ is a matching of $G$ if any two distinct edges of $M$ have no vertex in common.  Let $M= \{\{a_i,b_i\} \mid i = 1,\ldots,s\}$ be a nonempty matching of $G$. According to \cite{CV}, we say that $M$ is an {\it ordered matching} if:
\begin{enumerate}
\item $\{a_1,\ldots,a_s\}$ is a set of independent vertices,
\item $\{a_i,b_j\}\in E(G)$ implies $i \le  j$.
\end{enumerate}
\vskip 0.3cm
 
\begin{exm} In  the graph  $G = C_4$,  the subset $M = \{ \{a_1=1, b_1=4\},  \{2,3\}\}$ is a matching, but not an ordered matching, since the first property above would imply $a_2 = 3$ and $b_2 =2$. Then, $\{a_2,b_1\} = \{3,4\} \in E(G)$ and the second property above would  not hold.

\centerline{\setlength{\unitlength}{0.6cm}
 \begin{picture}(16,5)
\put(10,4){\line(2,0){2}}
\put(9.6,4){$1$}
\put(12,4){\line(1,-2){1}}
\put(12.1,4){$2$}
\put(13,2){\line(-2,-1){2}}
\put(13.1,2){$3$}
\put(11,1){\line(-2,1){2}}
\put(11,0.5){$4$}
\put(9,2){\line(1,2){1}}
\put(8.6,2){$5$}
\put(11,-0.6){$C_5$}
\put(1,4){\line(2,0){2}}
\put(0.6,4){$1$}
\put(3,4){\line(0,-2){2}}
\put(3.1,4){$2$}
\put(3,2){\line(-2,0){2}}
\put(3.1,1.5){$3$}
\put(1,2){\line(0,2){2}}
\put(0.6,1.5){$4$}
\put(2,0.5){$C_4$}
\end{picture}}
\vskip1cm  
 \noindent In the graph $G= C_5$,  by setting $a_1=1,\ b_1=2,\  a_2 = 4, \ b_2 = 3$, $M = \{\{1,2\}, \{3,4\} \}$ is  an ordered matching.
\end{exm}

\begin{defn}\label{ORDERED-MATCHING} The {\it ordered matching number} of $G$ is:
$$\nu_0(G) := \max \{|M| \mid M \subseteq E(G) \text{ is an ordered matching of }G\}.$$
\end{defn}

Then we have

\begin{thm} {\rm  \cite[Theorem 3.4]{HKTT}} \label{SPThm2} Let $G$ be a simple graph with $r$ vertices. Then,
$$\depth R/J(G)^{(n)} = r - \nu_0(G)-1 \text{ for all } n \ge  2\nu_0(G)-1.$$
In particular $\dstab^*(I) \le  2\nu_0(G)-1 < r$.
\end{thm}

The proof of this theorem is based on Takayama's Theorem \ref{Takay} and  Lemma  \ref{L02}. Note that this bound is sharp, see \cite[Proposition 3.6]{HKTT}.

In \cite{HeV} Herzog and Vladoiu describe  some classes of square-free monomial ideals $I$ with constant depth function, i.e. their $\dstab (I) = 1$.

\subsection{Depth of  powers }

Unlike the case of  integral closures of an arbitrary monomial ideal and symbolic powers of square-free monomial ideals, the behavior of depth function of a monomial ideal is very bad until it reaches the stability. This is first observed by Herzog and Hibi \cite{HH1}.  A more complicated picture is given in \cite{BHH}. Very recently Ha et al. obtain the following surprising result, which completely solves the problem of the  initial behavior of the depth function.

 \begin{thm}{\rm \cite[Theorem 6.7]{HHTT}} \label{HHTTThm} Let $f(n)$ be any convergent non-negative numerical function. Then there exists  a monomial ideal $I$ in $R=K[x_1,...,x_r]$ such that $f(n) = \depth (R/I^n)$ for all $n$.
\end{thm}

This maybe is a reason why Problem 1 for $\dstab(I)$ is much more difficult than for $\dstabIt(I)$. Example  3.1 in \cite{Hoa} constructed from Example 
\ref{Example1} shows that the bound (if exists) must be at least of the order $O(d(I)^{r-2})$.
However for the case of edge ideals there is a nice bound established by T. N. Trung \cite{NT2}. Recall that  a {\it leaf} in a graph $G$ is a vertex of degree one and a {\it leaf edge} is an edge incident with a leaf.  For an example, in Figure 3,  leafs are:   the vertex $4$ in $G_1$ and the vertex $2$ in $G_2$, and  edge leafs are:  the edge $\{1,4\}$ in $G_1$ and $\{ 2, 4\}$ in $G_2$. A connected graph is called a {\it tree} if it contains no cycles.  We use the symbols $\varepsilon (G)$ and $\varepsilon_0(G)$ to denote the number of edges and leaf edges of $G$, respectively.

 \begin{thm}{\rm \cite[Theorem  4.6]{NT2}}.   \label{NT2Thm1} Let  $G_1,...,G_s$ be all connected bipartite components  and $G_{s+1}, ..., G_{p}$   all connected non-bipartite  components of $G$. Let $2k_i$ be the maximal length of  cycles of $G_i,\ i\le s$ ($k_i = 1$  if $G_i$ has no cycle), and let $2k_i - 1$ be the maximal length of  odd cycles of $G_i,\ i> s$. 
 Then
 $$\dstabIt(I) \le r- \varepsilon_0(G) - \sum_{i=1}^p k_i +1 \le r - p.$$
 \end{thm}
 
 The proof is quite long and complicate.  First, the author studies connected graphs. From properties of $\ass(R/I(G)^n)$  it turns out that $\depth(R/I(G)^n) >0$ for all $n\ge 1$ provided $G$ is bipartite (see Theorem \ref{SVVThm}) and $\depth(R/I(G)^n) = 0$ for all $n\gg 0$  if $G$ is non-bipartite \cite[Corollary 3.4]{CMS}.  In the case of connected non-bipartite  graphs, the proof intensively uses the construction developed in \cite{CMS}. When  $G$ is a connected bipartite graph, thanks  to Theorem \ref{SVVThm},  $I(G)^n = I(G)^{(n)}$ for all $n\ge 1$. Hence, one can apply Lemma \ref{L02} to describe $\Da(I(G)^n)$. A key point in \cite{NT2} is to show that  $\depth(R/I(G)^n) = 1$ for all $n\gg 0$, see \cite[Lemma 3.1 and 3.3]{NT2}.   There the Takayama's Theorem \ref{Takay} is used only to show that  $\Da(I(G)^n))$ is disconnected, but Lemma \ref{L02} is very important. Some results in graph theory are also needed. 
Finally, the following result allows to  reduce the problem to connected components of $G$.
 
 \begin{thm} {\rm \cite[Theorem 4.4]{NT2}} Keep the notation in  Theorem \ref{NT2Thm1}. Then
 \begin{enumerate}
\item $\min \{\depth (R/ I(G)^n)|\ n\ge 1\} =s$.
\item $\dstab(I(G)) = \min\{n \ge 1 |\ \depth (R/I(G)^n) =s \}$.
\item $\dstab(I(G)) =\sum_{i=1}^p \dstab(R/I(G_i)) − p +1$.
\end{enumerate}
\end{thm}

Below are some other partial solutions to Problem 1.

\begin{enumerate} 
\item \cite[Corollary 3.4]{HH1}. The {square-free Veronese ideal} of degree $d$ in the variables $X_{i_1},...,X_{i_s},\ s\le r$, is the ideal of $R$ generated by all square-free monomials in $X_{i_1},...,X_{i_s}$ of degree $d$. Let $2\le d< n$ and let $I=I_{r,d}$ be the square-free Veronese ideal of degree $d$ in $r$ variables. Then 
$$\depth(R/I_{r,d}^n) = \max\{ 0, r- n(r-d) - 1\}.$$
In particular, $\dstab(I_{r,d}) \le \frac{r-1}{r-d}$.
\vskip0.3cm

\item \cite[Corollary 3.8]{HH1}.  Let $P$ be a finite partially ordered set (called {\it  poset} for short). A {\it poset  ideal} of  $P$ is a subset $I \subset P$ such that if $x\in I$, $y\in P$ and $y\le x$, then  $y\in P$. In particular, the empty set as well as $P$ itself is a poset ideal of $P$.  Write $\Jcal(P)$ for the finite  poset which consists of all poset ideals of  $P$, ordered by inclusion. 

Let $P = \{p_1,..., p_r\}$ be a finite poset and  $S= K[X_1,...,X_r,Y_1,...,Y_r]$. Consider the square-free  monomial ideal
$$H_P = (u_I := (\Pi_{p_i\in I}X_i)(\Pi_{p_i\not\in I}Y_i)|\ I\in \Jcal(P)) \subset  S.$$
Then $\dstab(H_P) = \rank(P) + 1 \le r$ and $\depth (S/H_P) > \depth (S/H_P^2)> \cdots > \depth (S/H_P^{\rank(P) + 1}) = r-1$, where $\rank(P)$ is the so-called rank of $P$. 
\vskip0.3cm

\item  \cite[Theorem 4.1]{HQ}. Let $I$ be a polymatroidal ideal. Then $\dstab(I) \le \ell(I)\le r$.
\vskip0.3cm

\item Let $\Hcal = (V, \Ecal)$ be a hypergraph. The incident mtrix $M= (a_{ij})$ of $\Hcal$ has $|V|$ rows and $|\Ecal|$ columns such that $a_{ij}= 1$ if $i \in E_j$ and $a_{ij} = 0$ otherwise. A hypergraph $\Hcal$ is said to be {\it unimodular} if  every square submatrix  of its incident matrix  has determinant equal to $0,1$ or $-1$.  Then \cite[Theorem 2.3 and Theorem 3.2]{HagT} state that $\depth (R/J(\Hcal)^n)$ is non-decreasing and $\dstab(J(\Hcal)) \le r$,  provided $\Hcal$ is an unimodular hypergraph.
\end{enumerate}

\subsection{Cohen-Macaulay property of powers} \label{CMness}

In this case, the depth gets its maximal value, and obtained results look nicest. Recall that a 
 {\it matroid} is a simplicial complex $\Delta$ with the following property: If $F, G \in \Delta$ and $|F| > |G|$, then there is $a\in F\setminus G$ such that $G\cup \{a\} \in \Delta $. 
 
 \begin{thm} {\rm  \cite[Theorem 1.1]{TT1}} \label{TTThm1} Let $\Delta$ be a simplicial complex. Then the following conditions are equivalent:
\begin{enumerate}
\item $R/I_{\Delta}^{(n)}$ is Cohen-Macaulay for every $n\ge  1$;
\item $R/I_{\Delta}^{(n)}$ is Cohen-Macaulay for some $n\ge  3$;
\item $\Delta$ is a matroid.
\end{enumerate}
\end{thm}

The equivalence $(1) \Leftrightarrow (3)$ are  independently proved by Minh and N.V. Trung \cite[Theorem 3.5]{MT2} and by Varbaro \cite[Theorem 2.1]{Var}. In the approach of \cite{MT2}, Lemma \ref{L02} plays an important role. It allows them to use tool from linear programming to show that the Cohen-Macaulayness of all symbolic powers characterizes matroid complexes. 

We say that $\Delta $ is a complete intersection if $I_\Delta $ is a complete intersection. This is equivalent to the property that no two minimal non-faces of $\Delta $ share a comon vertex. Since $R/I^n$ is Cohen-Macaulay if and only if $I^n = I^{(n)}$ and $R/I^{(n)}$ is Cohen-Macaulay, using Theorem \ref{TTThm1}, Terai and N. V. Trung can prove

\begin{thm} {\rm \cite[Theorem 1.2]{TT1})} \label{TTThm2}  Let $\Delta$ be a simplicial complex. Then the following conditions are equivalent:
\begin{enumerate}
\item $R/I_{\Delta}^n$ is Cohen-Macaulay for every $n\ge  1$;
\item $R/I_{\Delta}^n$ is Cohen-Macaulay for some $n\ge  3$;
\item $\Delta$ is a complete intersection. 
\end{enumerate}\end{thm}
The idea for the proof of $(2) \Rightarrow (3)$ in both theorems above  comes from the fact that matroid  and complete intersection complexes can be characterized by properties of their links. The main technical result of \cite{TT1}  shows that a complex $\Delta $ with
$\dim \Delta \ge 2$  is a matroid if and only if it is connected and locally a matroid. A similar result on complete intersection are proved in \cite[Theorem 1.5]{TY}. 

It is worth to mention   that the Cohen-Macaulay property of the second (ordinary or symbolic) power of a Stanley-Resiner ideal is completely different and is still not completely understood, see  \cite{HgT,  HgMT,  RTY, TT2, TY}.

The following result follows from Lemma \ref{DepthInc} and Theorem \ref{DStab}.

\begin{thm} {\rm \cite[Theorem 4.1]{HT3}}\label{CM} Let $I$ be a monomial ideal of $R$. The following conditions are equivalent
\begin{enumerate}
\item $R/\overline{I^n}$ is a Cohen-Macaulay ring for all $n\ge  1$,
\item  $R/\overline{I^n}$ is a Cohen-Macaulay ring for some $n\ge  r(r^2-1)r^{r/2}(r-1)^rd(I)^{(r-2)(r+1)}$,
\item $I$ is an equimultiple ideal of $R$. 
\end{enumerate}
\end{thm}

Following the idea of the proof of Theorem \ref{TTThm2} we also get a similar result for the integral closures.

\begin{thm} {\rm \cite[Theorem 4.7]{HT3}}\label{CMDelta}  Let $\Delta$ be a simplicial complex. Then the following conditions are equivalent:
\begin{enumerate}
\item $R/ \overline{I_{\Delta}^n}$ is Cohen-Macaulay for every $n\ge  1$;
\item $R/ \overline{I_{\Delta}^n}$ is Cohen-Macaulay for some $n\ge  3$;
\item $I_\Delta$ is a complete intersection;
\item $I_{\Delta}$ is an equimultipe ideal.
\end{enumerate}
\end{thm}

The property that the Cohen-Macaulayness  of $R/\overline{I^n}$ for some $n\ge  3$   forces that for all $n$ is very specific for  square-free monomial ideals. For an arbitrary monomial ideal, the picture is much more complicate, as shown by the following example.

\begin{exm}  \cite[Example 1]{HT3}  Let $d\ge  3$ and $I=(X^d,XY^{d-2}Z,Y^{d-1}Z) \subset R=K[X,Y,Z]$. Then 
\begin{enumerate}
\item $R/\overline{I^n}$ is Cohen-Macaulay for each $n=1,\ldots,d-1$;
\item $R/\overline{I^n}$ is not Cohen-Macaulay for any $n\ge  d$.
\end{enumerate}
Note that $\height(I) = 2$ and $\ell(I) = 3$ in this case.
\end{exm}

\subsection*{Acknowledgment} This work is partially supported by  the Project\\  VAST.HTQT.NHAT.01/16-18. I  would like to thank Prof. K. P. Shum for inviting me to be a keynote speaker at the Third  International Congress in Algebras and Combinatorics (ICAC 2017), Hong Kong, where I had chance to give this lecture to a broad audience.

\end{document}